\documentclass[11pt]{article}
\usepackage{amssymb}
\usepackage{amsfonts}
\usepackage{amscd}
\usepackage{epsfig}
\usepackage{amsmath,amstext,amsthm,amsbsy,amssymb}

\newcommand{\bc}{{\mathbb C}}
\newcommand{\br}{{\mathbb R}}
\newcommand{\bh}{{\mathbb H}}

\newcommand{\mtx}[4]

\begin{document}

\title{J{\o}rgensen's
inequality for quternionic hyperbolic space with elliptic elements
\author{WENSHENG CAO\thanks{Supported by NSFs of China (No.10801107,No.10671004 ),
NSF of Guangdong Province (No.8452902001000043) and Educational
Commission of Guangdong Province (No. LYM08097).} \ and\ \  HAIOU
TAN\thanks{Supported by NSF of Guangdong Province
(No.8152902001000004).}} }

\date{}
\maketitle
\bigskip
{\bf Abstract.} In this paper, we give an analogue of J{\o}rgensen's
inequality for non-elementary groups of isometries of quaternionic
hyperbolic space generated by two elements, one of which is
elliptic. As an application, we obtain an analogue of J{\o}rgensen's
inequality in 2-dimensional M\"obius group of the above case.
\smallskip

\medskip

\medskip


\newtheorem{thm}{Theorem}[section]
\newtheorem{defi}{Definition}[section]
\newtheorem{lem}{Lemma}[section]
\newtheorem{pro}{Proposition}[section]
\newtheorem{cor}{Corollary}[section]
\newtheorem{rem}{Remark}[section]

\section{Introduction}

 \quad\quad  J{\o}rgensen's  inequality \cite{jor} gives a necessary
condition for a non-elementary two generator subgroup of
$PSL(2,\bc)$ to be discrete.   Viewing  $PSL(2,\br)$  as the
isometry group of complex hyperbolic 1-space, $H_{\bc}^1$, one can
seek to generalize J{\o}rgensen's inequality to higher dimensional
complex hyperbolic isometries. There has been much research in this
area.
 Kamiya \cite{kam83,kam91} and Parker
\cite{par92,par97} gave generalizations of J{\o}rgensen's inequality
to the two generator subgroup of $PU(n,1)$ when one generator is
Heisenberg translation. By using the stable basin theorem, Basmajian
and Miner \cite{bm} generalized the J{\o}rgensen's  inequality to
the two generator subgroup of $PU(2, 1)$ when the generators are
loxodromic or boundary elliptic, and several other inequalities are
due to Jiang, Kamiya and Parker \cite{jkp} by using matrix method
other than the purely geometric method.   Jiang \cite{jp} and Kamiya
\cite{kp} gave generalizations of J{\o}rgensen's inequality to the
two generator subgroup of $PU(2,1)$ when one generator is Heisenberg
screw motion. The generalization also was done  in \cite{xj} for the
case when one generator is a regular elliptic element.

Following the research on complex hyperbolic space, I. Kim and J.
Parker opened up the study of quaternionic hyperbolic space in
\cite{kimp}.  They proved some basic facts about discreteness of two
generator subgroups, minimal volume of cusped quaternionic
manifolds, and laid some basic tools to study quaternionic
hyperbolic space.

It is naturally asked that theorems in complex hyperbolic space can
be generalized to quaternionic hyperbolic space. There is an attempt
in this area.  Kim \cite{kim} found  analogues in quaternionic
hyperbolic space of results in \cite{jkp,jp}.

The purpose of this paper is to show a condition for two generator
subgroup of $Sp(n,1)$ with an elliptic element to be not discrete.

In order to state our theorem, we recall some facts about elliptic
elements in $Sp(n,1)$. Every eigenvalue of elliptic element $g\in
Sp(n,1)$ has positive or negative type \cite{chen} and its
eigenvalues fall into $n$ similarity classes of positive type and
one similarity class of negative type. Let $\Lambda_i,\
i=1,\cdots,n$ be  its positive classes of eigenvalues and
$\Lambda_{n+1}$ be its negative class. Then any element in
$\Lambda_i$ has norm 1 and the fixed point set $F(g)$ of $g$ in
$H_{\bh}^n$, contains only one fixed point if $\Lambda_{n+1}\neq
\Lambda_i,\ i=1,\cdots,n$ and is a totally geodesic submanifold
which is equivalent to $H^m_{\bh}$ or $H^m_{\bc}$ (for some $m\leq
n)$ if $\Lambda_{n+1}$ coincides with exactly $m$ of the class
$\Lambda_i,\ i=1,\cdots,n$.  We call $g$  a {\it regular elliptic}
element if $F(g)$ contains only one point; Otherwise $g$ is called a
{\it boundary elliptic} element.  We mention here that the
definition of regular elliptic element is slightly different from
Goldman's \cite{gold} which requires its eigenvalues are distinct in
the setting of complex number. If $g\in Sp(n,1)$ is elliptic, then
$g$ conjugates to
$$diag(\lambda_1,\  \lambda_2,\  \cdots,\  \lambda_{n+1}),\eqno(1.1)$$
where $\lambda_i\in \Lambda_i,\ i=1,\cdots,n+1.$ We define
$$\delta(g)=\max\{|\lambda_i-\lambda_{n+1}|^2: \ i=1,\cdots,n\}.\eqno(1.2)$$ Since
similarity classes $\Lambda_i,\  i=1,\cdots,n+1$ of $g$ are
invariant under conjugation, we have the following proposition.
\begin{pro}
If $g\in Sp(n,1)$  is  elliptic, then $\delta(g)$ is invariant under
conjugation.
\end{pro}

We now deduce the formula of  $\delta(g)$ defined by (1.2) for an
elliptic element $g\in Sp(n,1).$

We use \cite{zhang} as a reference for the properties of quaternions
and matrices of quaternions, also see  in Section 2 for an
abbreviated description. Since each similarity class $\Lambda_i$ has
a unique complex number with nonnegative imaginary part. Let
$$e^{{\bf i}\theta_i}\in \Lambda_i, \ i=1,\cdots,n+1$$ be such
complex numbers, that is, $0\leq \theta_i\leq \pi$ for each $i$. Let
$$\delta_{i,n+1}=\max\{|\lambda_i-\lambda_{n+1}|^2:\  \lambda_i\in \Lambda_i,\ \lambda_{n+1}\in \Lambda_{n+1}\},i=1,\cdots,n.$$
Then
$$\delta_{i,n+1}=\max_{u,w\in \bh}\{|ue^{{\bf i}\theta_i}u^{-1}-we^{{\bf
i}\theta_{n+1}}w^{-1}|^2\}=\max_{|w|=1}\{|e^{{\bf
i}\theta_i}-we^{{\bf i}\theta_{n+1}}w^{-1}|^2\}.$$  Let
$$T(w)=e^{{\bf i}\theta_i}\overline{w^{-1}} e^{{-\bf i}\theta_{n+1}}\bar{w}+we^{{\bf
i}\theta_{n+1}}w^{-1} e^{{-\bf i}\theta_{i}}.$$  Then $|e^{{\bf
i}\theta_i}-we^{{\bf i}\theta_{n+1}}w^{-1}|^2=2-T(w)$ and
$$\delta_{i,n+1}=2-\min_{|w|=1} T(w).$$
We take the complex representation of  the quaternion $w$ with norm
1  as
$$w =w_1+w_2{\bf
j},\ \ \mbox{where}\ \ w_1,w_2\in \bc.$$ Then
$w^{-1}=\bar{w}=\overline{w_1}-w_2{\bf j}$. Note that
$$z{\bf j}={\bf j}\bar{z}\ \mbox{for}\ z\in \bc.$$  By direct computation we
have
$$T(w)=2(|w_1|^2cos(\theta_i-\theta_{n+1})+|w_2|^2cos(\theta_i+\theta_{n+1})).$$
Hence\begin{eqnarray*} \min_{|w|=1}T(w)=\Big\{\begin{array}{l}
2\cos(\theta_i+\theta_{n+1})\, \ \ \mbox{if}\ \
\cos(\theta_i-\theta_{n+1})\geq \cos(\theta_i+\theta_{n+1})
\\ 2\cos(\theta_i-\theta_{n+1})\,\ \  \mbox{if}\ \
\cos(\theta_i-\theta_{n+1})< \cos(\theta_i+\theta_{n+1}).
\end{array}
\end{eqnarray*}
Therefore
$$\delta(g)=\max\{4\sin^2\frac{\theta_i\pm\theta_{n+1}}{2}:\ i=1,\cdots,n\}.\eqno(1.3)$$

Our main theorem is
\begin{thm}\label{main} Let $g $ and $h$ be elements of $Sp(n,1).$ Suppose that $g$
is an elliptic element with fixed point set $F(g)\in H_{\bh}^n.$  If
$$\inf_{q\in F(g)} \cosh ^{2} {{\rho(q, h(q))} \over 2}\delta(g) < 1,\eqno(1.4)$$ then
the group $<g, h>$ generated by $g $ and $h$ is either elementary or
not discrete.
\end{thm}

Applying Theorem \ref{main} to the subgroup $PU(2,1)$ of $Sp(2,1)$
with $g$ being regular elliptic, we obtain
\begin{cor}(cf.\cite[Theorem 3.4]{xj})
Let g and h be elements of $PU(2,1)$ so that g is a regular elliptic
element with unique fixed point $q$.  If
$$ cosh^2\frac{\rho(q, h(q))}{2}\delta(g)<1,$$
then the group $\langle g,h \rangle$  is either elementary or not
discrete.
\end{cor}

Let $$g=\left(
          \begin{array}{cc}
            e^{{\bf i}\theta} & 0 \\
            0 &  e^{-{\bf i}\theta}\\
          \end{array}
        \right), \ \  h=\left(
                          \begin{array}{cc}
                            a & b \\
                            c & d \\
                          \end{array}
                        \right)\in SL(2,\bc).\eqno(1.5)$$

  As an application,  by
embedding $SL(2,\bc)$ to $Sp(1,1)$ in Section 4,  we obtain

\begin{thm}\label{mob}
Let g and h be elements of $SL(2,\bc)$ given by (1.5).  If $ \langle
g, h \rangle$ is discrete and non-elementary. Then
$$\inf_{|t|<1, t\in \bc}4f(t)\sin^2\theta\geq 1,  \eqno(1.6)$$
where
$$f(t)=\frac{|1-t^2|^2(|a|^2+|d|^2)+|1+t|^4|c|^2+|1-t|^4|b|^2+2(\bar{t}-t)^2}{4(1-|t|^2)^2}+\frac{1}{2}. \eqno(1.7)$$
\end{thm}

Choosing  $t=0$ in the above theorem, we obtain

\begin{cor}\label{mob1}
Let g and h be elements of $SL(2,\bc)$ given by (1.5).  If $ \langle
g, h \rangle$ is discrete and non-elementary. Then
$$\sin^2\theta(||h||^2+2)\geq 1,  \eqno(1.8)$$
where $||h||^2=|a|^2+|b|^2+|c|^2+|d|^2$.
\end{cor}

\begin{rem}
The J{\o}rgensen's inequality \cite{jor} gives that
$$|tr(g)^2-4|+|tr(ghg^{-1}h^{-1})-2|=4\sin^2\theta(1+|bc|)\geq 1$$
in the above case.

  Let $h=\left(  \begin{array}{cc}
                            \frac{-3}{2} & 2{\bf i} \\
                            2{\bf i} & 2 \\
                          \end{array}
                        \right)$. Then
                        $16\frac{1}{4}=(||h||^2+2)<4(1+|bc|)=20$,
                        which implies  that Corollary \ref{mob1} is better
                        than the J{\o}rgensen's inequality in such a case.   While the
                        case  $h=\left(  \begin{array}{cc}
                            1 & \sqrt{2} \\
                            \sqrt{2} & 3 \\
                          \end{array}
                        \right)$ having $$f(t)={\frac{7+14t_1^2t_2^2+2t_1^2+10t_2^2+7(t_1^4+t_2^4)}{2(1-t_1^2-t_2^2)^2}+\frac{1}{2}}>4>1+|bc|=3$$  implies Theorem \ref{mob} is weaker
                        than the J{\o}rgensen's inequality in such a case.

                        By the above  comparison of   Theorem 1.2 with the J{\o}rgensen's inequality, we show that neither theorem is
a consequence of the other one.


\end{rem}

\section{Preliminaries}

\quad\quad  In this section, we give some necessary background
materials of quaternionic  hyperbolic geometry. More details can be
found in \cite{chen,gold,kimp}.

Let $\bh$ denote the division ring of real quaternions. Elements of
$\bh$ have the form $z=z_1+z_2{\bf i}+z_3{\bf j}+z_4{\bf k}\in \bh$
where $z_i\in \br$ and
$$
{\bf i}^2 = {\bf j}^2 = {\bf k}^2 = {\bf i}{\bf j}{\bf k} = -1.
$$
Let $\overline{z}=z_1-z_2{\bf i}-z_3{\bf j}-z_4{\bf k}$ be the {\sl
conjugate} of $z$, and
$$
|z|= \sqrt{\overline{z}z}=\sqrt{z_1^2+z_2^2+z_3^2+z_4^2}
$$
be the {\sl modulus} of $z$. We define $\Re(z)=(z+\overline{z})/2$
to be the {\sl real part} of $z$, and $\Im(z)=(z-\overline{z})/2$ to
be the {\sl imaginary part} of $z$. Also
$z^{-1}=\overline{z}|z|^{-2}$ is the {\sl inverse} of $z$.  Observe
that $\Re(wzw^{-1})=\Re(z)$ and $|wzw^{-1}|=|z|$ for all $z$ and $w$
in $\bh$. Two quaternions $z$ and $w$ are {\sl similar} if there
exists nonzero $q \in \bh$ such that $z=q w q^{-1}$. The similarity
class of $z$ is the set $\bigl\{ q z q^{-1} : q \in \bh-\{ 0\}
\bigr\}$.

Let $\bh^{n,1}$ be the vector space of dimension n+1 over $\bh$ with
the unitary structure defined by the Hermitian form
$$
\langle{\bf z},\,{\bf w}\rangle={\bf w}^*J{\bf z}=
\overline{w_1}z_1+\cdots+\overline{w_n}z_n-\overline{w_{n+1}}z_{n+1},
$$
where ${\bf z}$ and ${\bf w}$ are the column vectors in $V$ with
entries $(z_1,\cdots,z_{n+1})$ and $(w_1,\cdots,w_{n+1})$
respectively, $\cdot^*$ denotes the conjugate transpose and $J$ is
the Hermitian matrix
$$J=\left(
      \begin{array}{cc}
        I_n & 0 \\
         0 & -1 \\
      \end{array}
    \right).$$
We define a {\sl unitary transformation} $g$ to be an automorphism
of $\bh^{n,1}$, that is, a linear bijection such that $\langle
g({\bf z}),\,g({\bf w})\rangle=\langle{\bf z},\,{\bf w}\rangle$ for
all ${\bf z}$ and ${\bf w}$ in $V$. We denote the group of all
unitary transformations by $Sp(n,1)$.

Following Section 2 of \cite{chen}, let
\begin{eqnarray*}
V_0 & = & \Bigl\{{\bf z} \in  V-\{0\}:
\langle{\bf z},\,{\bf z}\rangle=0\Bigr\} \\
V_{-} &  = & \Bigl\{{\bf z} \in V:\langle{\bf z},\,{\bf
z}\rangle<0\Bigr\}.
\end{eqnarray*}
It is obvious that $V_0$ and $V_{-}$ are invariant under $ Sp(n,1)$.
We define $V^s$ to be $V^s=V_{-}\cup  V_0$. Let $P:V^s\to
P(V^s)\subset \bh^{n}$ be the projection map defined by
$$
P(z_1,\cdots,z_n,  z_{n+1})^t=(z_1z_{n+1}^{-1},\cdots,z_n
z_{n+1}^{-1})^t,
$$
where $\cdot^t$ denotes the  transpose.

 We define $H_{\bh}^n=P(V_-)$
and $\partial H_{\bh}^n=P(V_0)$.
 The Bergman metric on $H_{\bh}^n$ is given by the distance formula
 $$\cosh^2\frac{\rho(z,w)}{2}=\frac{\langle{\bf z},\,{\bf
w}\rangle \langle{\bf w},\,{\bf z}\rangle}{\langle{\bf z},\,{\bf
z}\rangle \langle{\bf w},\,{\bf w}\rangle},\ \ \mbox{where}\ \ {\bf
z}\in P^{-1}(z),{\bf w}\in P^{-1}(w). \eqno(2.1)$$
 The holomorphic isometry group of $H_{\bh}^n$ with respect to the
Bergman metric is the projective unitary group $PSp(n, 1)$ and acts
on $P(\bh^{n,1})$ by matrix multiplication.

If $g\in Sp(n,1)$, by definition, $g$ preserves the Hermitian form.
Hence
$$
{\bf w}^*J{\bf z}=\langle{\bf z},\,{\bf w}\rangle= \langle g{\bf
z},\,g{\bf w}\rangle ={\bf w}^* g^*J g{\bf z}
$$
for all ${\bf z}$ and ${\bf w}$ in $V$. Letting ${\bf z}$ and ${\bf
w}$ vary over a basis for $V$ we see that $J= g^*J g$. From this we
find $ g^{-1}=J^{-1} g^*J$. That is:
$$
g^{-1}=\left(
  \begin{array}{cc}
     A^*& -\beta^* \\
    -\alpha^* & \overline{a_{n+1,n+1}}\\
    \end{array}
\right)\ \  \mbox{for} \ \ g=(a_{ij})_{i,j=1,\cdots,n+1}=\left(
  \begin{array}{cc}
     A& \alpha \\
    \beta& a_{n+1,n+1}\\
    \end{array}
\right).$$ Using the identities $gg^{-1}=g^{-1}g=I$ we obtain:

$$ AA^*-\alpha \alpha^*=I_n,\ \  -A\beta^*+\alpha
\overline{a_{n+1,n+1}}= 0,\ \
-|\beta|^2+|a_{n+1,n+1}|^2=1;\eqno(2.2)$$
$$
A^*A-\beta^*\beta=I_n,\ \ A^*\alpha-\beta^*a_{n+1,n+1}= 0,\ \
-|\alpha|^2+|a_{n+1,n+1}|^2=1.\eqno(2.3)
$$

For a non-trivial element $g$ of $Sp(n,1)$, we say that $g$ is {\it
parabolic} if it has exactly one fixed point and this lies on
$\partial H_{\bh}^n$,  $g$ is {\it loxodromic} if it has exactly two
fixed points and they lie on $\partial H_{\bh}^n$ and $g$ is {\it
elliptic} if it has a fixed point in $H_{\bh}^n$. In particular, if
$g$ has fixed point $q_0=(0, \cdots, 0)^t\in H_{\bh}^n$, then $g$
has the form
$$g=diag(A,a),
$$
where $A\in U(n;\bh)$ and $a\in U(1;\bh)$.

 A  subgroup $G$ of
$Sp(n,1)$ is called {\it non-elementary} if it contains two
non-elliptic elements of infinite order with distinct fixed points;
Otherwise $G$ is called {\it elementary}.

As in complex hyperbolic n-space, we have the following proposition
classifying  elementary subgroups of $Sp(n, 1)$.
\begin{pro}(cf.\cite[Lemma 2.4]{caobull} \label{element}
(i)\quad If  $G$ contains a parabolic element but no loxodromic
element, then $G$ is elementary if and only if it fixes a point in
$\partial H_{\bh}^n$;

(ii)\quad If  $G$ contains a loxodromic element, then $G$ is
elementary if and only if it fixes a point in $\partial H_{\bc}^n$
or a point-pair $\{x, y\}\subset \partial H_{\bh}^n$;

(iii)\quad $G$ is purely elliptic, i.e., each non-trivial element of
$G$ is elliptic, then  $G$ is elementary and fixes a point in
$\overline{H_{\bh}^n}$.
\end{pro}

\section{Proof of Theorem \ref{main}}

{\bf Proof of Theorem \ref{main}.} \ \   For  any fixed point $q\in
fix(g)$, since  $cosh^2\frac{\rho(q, h(q))}{2}\delta(g)$ is
invariant under conjugation, we may assume that g is of the form
(1.1) having fixed point  $q=(0, \cdots, 0)^t\in H_{\bh}^n$ and
$$h=(a_{ij})_{i,j=1,\cdots,n+1}=\left(
  \begin{array}{cc}
     A& \alpha \\
    \beta& a_{n+1,n+1}\\
    \end{array}
\right).$$
Then
$$cosh^2\frac{\rho(q, h(q))}{2}=|a_{n+1,n+1}|^2,\ \ \delta(g)=\max\{|\lambda_i-\lambda_{n+1}|^2: \ i=1,\cdots,n\} .$$

In what follows, we will show that if  $$ |a_{n+1,n+1}|^2 \delta(g)<
1, \eqno(3.1)$$ then  the group $\langle g,h \rangle$  is either
elementary or not discrete.

 Let $h_0 = h$  and  $h_{k+1} = h_kgh_k^{-1}$.  We
write
$$h_k=(a^{(k)}_{ij})_{i,j=1,\cdots,n+1}=\left(
  \begin{array}{cc}
     A^{(k)}& \alpha^{(k)} \\
    \beta ^{(k)}& a^{(k)}_{n+1,n+1}\\
    \end{array}
\right).$$
 Then
\begin{eqnarray*}
h_{k+1}&=& \left(
  \begin{array}{cc}
     A^{(k+1)}& \alpha^{(k+1)} \\
    \beta ^{(k+1)}& a^{(k+1)}_{n+1,n+1}\\
    \end{array}
\right)\\
 &=& \left(
  \begin{array}{cc}
     A^{(k)}& \alpha^{(k)} \\
    \beta ^{(k)}& a^{(k)}_{n+1,n+1}\\
    \end{array}
\right)\left(
       \begin{array}{cc}
         L & 0 \\
         0 & \lambda_{n+1} \\
            \end{array}
     \right)\left(
  \begin{array}{cc}
     (A^{(k)})^* & -(\beta^{(k)})^* \\
    -(\alpha^{(k)})^*& \overline{a^{(k)}_{n+1,n+1}}\\
    \end{array}
\right),
\end{eqnarray*}
where $L=diag(\lambda_1,\  \lambda_2,\  \cdots,\  \lambda_{n}).$

\noindent Therefore $$ a^{(k+1)}_{n+1,n+1} =
a^{(k)}_{n+1,n+1}\lambda_{n+1}\overline{a^{(k)}_{n+1,n+1}}-\beta
^{(k)}L(\beta^{(k)})^* \eqno (3.2)$$ and
\begin{eqnarray*}
|a^{(k+1)}_{n+1,n+1}|^2
&=&(a^{(k)}_{n+1,n+1}\lambda_{n+1}\overline{a^{(k)}_{n+1,n+1}}-\beta
^{(k)}L(\beta^{(k)})^*)(a^{(k)}_{n+1,n+1}\overline{\lambda_{n+1}}\overline{a^{(k)}_{n+1,n+1}}-\beta
^{(k)}L^*(\beta^{(k)})^*)\\
&=& |a^{(k)}_{n+1,n+1}|^4+ \beta ^{(k)}L(\beta^{(k)})^*\beta
^{(k)}L^*(\beta^{(k)})^*  \hspace{70mm}(3.3)\\
&&-a^{(k)}_{n+1,n+1}\lambda_{n+1}\overline{a^{(k)}_{n+1,n+1}}\beta
^{(k)}L^*(\beta^{(k)})^*-\beta
^{(k)}L(\beta^{(k)})^*a^{(k)}_{n+1,n+1}\overline{\lambda_{n+1}}\overline{a^{(k)}_{n+1,n+1}}.
\end{eqnarray*}

If there exists  some $k$ such that $\beta ^{(k)}=0$, then by (2.2)
and (2.3) we have $$|a^{(k)}_{n+1,n+1}|=1 \ \ \mbox{and}\ \ \alpha
^{(k)}=0,$$ which implies that $q$ is a fixed point of
$h_k=h_{k-1}gh_{k-1}^{-1}$.  We deduce that $q$ is a fixed point of
$h_{k-1}$  and  $|a_{n+1,n+1}|=1$ by induction. Thus  $q$ is a fixed
point of $h$, which implies that $\langle g,h \rangle$ is
elementary.

In what follows, we may assume that $\beta ^{(k)} \neq 0$.

We first consider the case when all the elements of $\beta ^{(k)}$
are nonzero, that is, $$a^{(k)}_{n+1,i}\neq 0, i=1,\cdots, n.$$

 In this case,
noting that $$\beta ^{(k)}L(\beta^{(k)})^*\beta
^{(k)}L^*(\beta^{(k)})^*\leq  |\beta ^{(k)}|^4,$$ we have
\begin{eqnarray*}
|a^{(k+1)}_{n+1,n+1}|^2 \leq && |a^{(k)}_{n+1,n+1}|^4+ |\beta
^{(k)}|^4-a^{(k)}_{n+1,n+1}\lambda_{n+1}\overline{a^{(k)}_{n+1,n+1}}(\sum_{i=1}^{n}a^{(k)}_{n+1,i}\overline{\lambda_{i}}\overline{a^{(k)}_{n+1,i}}
)\\
&&-(\sum_{i=1}^{n}a^{(k)}_{n+1,i}\lambda_{i}\overline{a^{(k)}_{n+1,i}})a^{(k)}_{n+1,n+1}\overline{\lambda_{n+1}}\overline{a^{(k)}_{n+1,n+1}}.\hspace{60mm}(3.4)
\end{eqnarray*}
 Let
$$u_i=\overline{a^{(k)}_{n+1,i}} \ ^{-1}
\lambda_i \overline{a^{(k)}_{n+1,i}}, \ i=1,\cdots,n+1. $$ Then
$u_i\in \Lambda_i, \ i=1,\cdots,n+1$ and
$$u_i\overline{u_{n+1}}+u_{n+1}\overline{u_i}=2-|u_i-u_{n+1}|^2.$$

We can rewrite (3.4) as $$ |a^{(k+1)}_{n+1,n+1}|^2 \leq
|a^{(k)}_{n+1,n+1}|^4+ |\beta
^{(k)}|^4-\sum_{i=1}^{n}|a^{(k)}_{n+1,n+1}|^2|a^{(k)}_{n+1,i}|^2(2-|u_i-u_{n+1}|^2).\eqno(3.5)
$$
 Noting that $-|\beta^{(k)}|^2+|a^{(k)}_{n+1,n+1}|^2=1$, by (3.5) we
 have
$$ |a^{(k+1)}_{n+1,n+1}|^2-1\leq
|a^{(k)}_{n+1,n+1}|^2\sum_{i=1}^{n}|a^{(k)}_{n+1,i}|^2|u_i-u_{n+1}|^2.
$$
Therefore $$|a^{(k+1)}_{n+1,n+1}|^2-1 \leq
(|a^{(k)}_{n+1,n+1}|^2-1)|a^{(k)}_{n+1,n+1}|^2\delta(g).\eqno(3.6)$$

We remark that the case $\beta ^{(k)} \neq 0$ with  some
$a^{(k)}_{n+1,t}=0$ for $t\in \{1,\cdots, n\}$   also leads to the
above inequality.

Noting (3.1),  we  obtain by induction
$$|a^{(k+1)}_{n+1,n+1}|<|a^{(k)}_{n+1,n+1}|$$ and
$$|a^{(k+1)}_{n+1,n+1}|^2-1 <
(|a_{n+1,n+1}|^2-1)(|a_{n+1,n+1}|^2\delta(g))^{k+1}.\eqno(3.7)$$
Thus $|a^{(k)}_{n+1,n+1}|\to 1$ and $\{h_k\}$ is a sequence with
distinct elements.  By (2.2) and (2.3), we have
$$\beta ^{(k)}\to 0,\ \ \alpha^{(k)}\to 0$$
and $$A^{(k)}(A^{(k)})^*\to I_n.$$ By passing to its subsequence, we
may assume
 $$A^{(k_t)}\to A_{\infty},\ \  a^{(k_t)}_{n+1,n+1}\to a_{\infty}.$$
 Thus $h_{k+1}$  converges to
 $$h_{\infty}= \left(
                 \begin{array}{cc}
                   A_{\infty} & 0 \\
                  0 & a_{\infty} \\
                  \end{array}
               \right)\in Sp(n,1),
$$
which implies that $\langle h, g\rangle$ is not discrete.  This
concludes the  proof.

\section{Proof of Theorem \ref{mob}}
{\bf Proof of Theorem \ref{mob}.} \ \ As in \cite{cpw,chen}, we can
regard $Sp(1,1)$ as the isometries of hyperbolic 4-space $H^4$,
whose model is the unit ball in the quaternions $\bh$.  $SL(2,\bc)$,
the isometries of hyperbolic 3-space $H^3$, can be embedded as a
subgroup of $Sp(1,1)$ as following:

$$f\in SL(2,\bc)\hookrightarrow TfT^{-1}\in Sp(1,1),\eqno(4.1)$$
where $$T=\frac{1}{\sqrt{2}}\left(
            \begin{array}{cc}
              1 & -{\bf j} \\
              -{\bf j} & 1 \\
            \end{array}
          \right).$$

Let $g$ and $h$ be stated as (1.5) and $\hat{g}=TgT^{-1},\,
\hat{h}=ThT^{-1}$.   Then
$$\hat{g}=\left(
          \begin{array}{cc}
            e^{i\theta} & 0 \\
            0 &  e^{-i\theta}\\
          \end{array}
        \right), \ \  \hat{h}=\frac{1}{2}\left(
                                                            \begin{array}{cc}
                                                              1 & -{\bf j} \\
                                                             -{\bf j} & 1 \\
                                                            \end{array}
                                                          \right)\left(
                                                                   \begin{array}{cc}
                                                                     a & b \\
                                                                     c & d \\
                                                                   \end{array}
                                                                 \right)
                                                          \left(
                                                            \begin{array}{cc}
                                                              1 & {\bf j} \\
                                                              {\bf j} & 1 \\
                                                            \end{array}
                                                          \right)\in Sp(1,1).\eqno(4.3)$$
In fact, $\hat{g}, \hat{h}\in Sp(1,1)$  can be verified directly by
lemma 1.1 in \cite{cpw}. Applying the  formula (1.3) to $\hat{g}$ in
which $\theta_1=\theta_2=\theta$, we have
$$\delta(\hat{g})=4\sin^2\theta.$$
It is easy to know that the fixed point set of $\hat{g}$ is $\{t{\bf
j}:t\in \bc \, \mbox{with}\,\, |t|<1\}.$

Let ${\bf z}=\left(
         \begin{array}{c}
           t{\bf j} \\
          1\\
         \end{array}
       \right)$ and ${\bf w}=\hat{h}{\bf z}$. Then
$$\langle{\bf w},\,{\bf z}\rangle=\frac{1}{2}(-t{\bf j},1)\left(
                                                            \begin{array}{cc}
                                                              1 & -{\bf j} \\
                                                              {\bf j} & -1 \\
                                                            \end{array}
                                                          \right)\left(
                                                                   \begin{array}{cc}
                                                                     a & b \\
                                                                     c & d \\
                                                                   \end{array}
                                                                 \right)
                                                          \left(
                                                            \begin{array}{cc}
                                                              1 & {\bf j} \\
                                                              {\bf j} & 1 \\
                                                            \end{array}
                                                          \right)\left(
         \begin{array}{c}
           t{\bf j} \\
          1\\
         \end{array}
       \right)$$
and
$$|\langle{\bf w},\,{\bf z}\rangle|^2=\frac{1}{4}((1-t){\bf j},-1-t)\left(
                                                                   \begin{array}{cc}
                                                                     a & b \\
                                                                     c & d \\
                                                                   \end{array}
                                                                 \right)
                                                          \left(
                                                            \begin{array}{cc}
                                                              |t|^2+1+\bar{t}+t & (1+t-\bar{t}- |t|^2){\bf j} \\
                                                              ( -1-t+\bar{t}+|t|^2){\bf j} &  |t|^2+1-\bar{t}-t\\
                                                            \end{array}
                                                          \right)\left(
                                                                   \begin{array}{cc}
                                                                     \bar{a} & \bar{c} \\
                                                                    \bar{b} & \bar{d} \\
                                                                   \end{array}
                                                                 \right)\left(
         \begin{array}{c}
           (t-1){\bf j} \\
         -1 -\bar{t}\\
         \end{array}
       \right).$$
Since $\langle{\bf z},\,{\bf z}\rangle=\langle{\bf w},\,{\bf
w}\rangle=|t|^2-1$, by direct computation, we have

 $$\cosh^2\frac{\rho(t{\bf j},\hat{h}(t{\bf j}))}{2}=\frac{\langle{\bf z},\,{\bf
w}\rangle \langle{\bf w},\,{\bf z}\rangle}{\langle{\bf z},\,{\bf
z}\rangle \langle{\bf w},\,{\bf w}\rangle}=f(t),$$ where $f(t)$ be
stated as in (1.7).

Applying  Theorem \ref{main} to  $\hat{g},\hat{h}\in  Sp(1,1)$ given
by (4.3), we conclude the proof of Theorem 1.2.

WENSHENG CAO and  HAIOU TAN,  Department of Mathematics and Physics,
  Wuyi University,Jiangmen  529020,  China\\
  \noindent e-mail: wenscao@yahoo.com.cn

\noindent e-mail: hotan@wyu.cn
\end{document}